\documentclass{birkmult}
\usepackage{amsmath, amssymb,color,ulem}
\usepackage{graphicx}

\newtheorem{theorem}{Theorem}

\newtheorem{lemma}[theorem]{Lemma}
\newtheorem{proposition}[theorem]{Proposition}

\newcommand{\ud}{\mathrm{d}}

\newcommand{\A}{\mathcal{A}}
\newcommand{\B}{\mathcal{B}}

\newcommand{\I}{\mathcal{I}}
\newcommand{\X}{\mathcal{X}}
\newcommand{\T}{\mathcal{T}}

\newcommand{\RR}{\mathbb{R}}

\renewcommand{\theequation}{\thesection.\arabic{equation}}

\newcommand{\sector}[1]{
  \refstepcounter{section}
  \setcounter{equation}{0}
  \setcounter{theorem}{0}
  \mbox{\ }\smallskip
  {\begin{center}
    \bf \thesection.\ {#1}
  \end{center}}}

\newcommand{\eofproof}{\hfill$\Box$ } 
\newcommand{\Pf}[1]{{\noindent \bf Proof{#1}\ }}

\begin{document}

% \title{On a size-structured two-phase population model with infinite states-at-birth\thanks{{\em 2000 MSC:} 92D25, 47D06, 35B35}}
% \author{J\'{o}zsef Z. Farkas\thanks{E-mail: jzf@maths.stir.ac.uk}\ {\normalsize $^1$}\mbox{\ }
% and Peter Hinow\thanks{E-mail: hinow@ima.umn.edu}\ {\normalsize $^2$} \mbox{\
% }
%     \mbox{\ }\\
%     \mbox{\ }\\
%     {\normalsize $^1$\,Department of Computing Science and Mathematics,
%      University of Stirling,}\\
%     {\normalsize  Stirling, FK9 4LA, UK}  \\
% {\normalsize $^2$\,Institute for Mathematics and its Applications, University of Minnesota,}\\
%     {\normalsize  114 Lind Hall, Minneapolis, MN 55455-0134, USA} }
\title[A size-structured two-phase population model]{On a size-structured two-phase population mo\-del with infinite states-at-birth}
\author[J\'{o}zsef Z. Farkas]{J\'{o}zsef Z. Farkas}
\address{Department of Computing Science and Mathematics, University of Stirling, Stirling, FK9 4LA, United Kingdom }
\email{jzf@maths.stir.ac.uk}

\author{Peter Hinow}
\address{Department of Mathematical Sciences, University of
Wisconsin -  Milwaukee,  P. O. Box 413, Milwaukee, WI 53201-0413,
USA}
\email{hinow@uwm.edu}

\date{\today}

\begin{abstract}
In this work we introduce and analyze a linear size-structured population model with infinite states-at-birth. 
We model the dynamics of a population in which individuals have two distinct life-stages: an ``active'' phase when individuals grow, 
reproduce and die and a second ``resting'' phase when individuals only grow. Transition between these two phases depends on  individuals' size. 
First we show that the problem is governed by a positive quasicontractive semigroup on the biologically relevant state space. 
Then we investigate, in the framework of the spectral theory of linear operators, the asymptotic behavior of solutions of the model. 
We prove that the associated semigroup has, under biologically plausible assumptions, the property of asynchronous exponential growth.
%{\em Keywords:} Size-structured populations; positivity; quasicontractive semigroups; spectral methods; asynchronous exponential growth.
\end{abstract}
\maketitle

\sector{Introduction}

Interest in understanding the dynamics of biological populations is old. Classical, ordinary differential equation models assume homogeneity 
of individuals within population classes, and involve equations for total population sizes. However, individuals in biological populations differ 
in their physiological characteristics. Vital rates, such as those of birth, death and development, vary amongst individuals. 
Therefore physiologically structured partial differential equation models are often more useful to understand the dynamics of biological populations. 
We refer the interested reader to the monographs \cite{I,MD,WEB} for basic concepts and results in the theory of structured populations. 

In this paper we study the following linear size-structured model
\begin{align}
u_{1,t}(s,t)+\left(\gamma_1(s)u_1(s,t)\right)_s & = -\mu(s)u_1(s,t)+\int_0^m\beta(s,y)u_1(y,t)\,\ud y  \nonumber \\ 
&\quad -c_1(s)u_1(s,t)+c_2(s)u_2(s,t),\label{model1}\\
u_{2,t}(s,t)+\left(\gamma_2(s)u_2(s,t)\right)_s & =c_1(s)u_1(s,t)-c_2(s)u_2(s,t).\label{model2}
\end{align}
Equations \eqref{model1}-\eqref{model2} are equipped with the following boundary and initial conditions
\begin{equation}\label{bd_in}
\begin{aligned}
& \gamma_1(0)u_1(0,t)=0,\quad u_1(s,0)=u^{*}_{1}(s),\\
& \gamma_2(0)u_2(0,t)=0,\quad u_2(s,0)=u^{*}_{2}(s).
\end{aligned}
\end{equation}
Individuals may experience two different stages in their life that we call ``active'' and ``resting''. 
The density of individuals in the active stage of size $s\in[0,m]$ at time $t$ is denoted by $u_1(s,t)$, 
while the density of individuals in the resting stage of size $s$ at time $t$ is denoted by $u_2(s,t)$. 
The maximal size an individual may reach is denoted by $m$. Individuals grow in both classes at the size-dependent 
growth rates $\gamma_1$ and $\gamma_2$, respectively. In the active stage individuals experience size-dependent mortality denoted by $\mu$. 
Further, only individuals in the active stage reproduce, this is expressed via the recruitment term on the right hand side of equation \eqref{model1}. 
In particular, the function $\beta(s,y)$ gives the rate at which an individual of size $y$ produces offspring of the size $s$.  
The transition between the two life-stages is captured by the size-dependent rates $c_1$ and $c_2$. We make the following assumptions on the model parameters 
\begin{equation}\label{assumptions}
\begin{aligned}
\mu,\, c_1,\, c_2 \in L_+^\infty([0,m]),&\quad  0\le\beta \in C([0,m]^2), \quad \beta\not\equiv 0,  \quad 0< \gamma_1,\,\gamma_2\in C^1([0,m]).
\end{aligned}
\end{equation}
Note that our assumptions in \eqref{assumptions} are tailored toward the mathematical analysis of model \eqref{model1}-\eqref{bd_in}. 
In case of a specific population, one can make additional assumptions on the model ingredients, such as $\beta(s,y)=0$ whenever $s>y$, 
that is, individuals can only produce offspring of smaller size. For later use, let 
\begin{equation*}
B=||\beta||_\infty\quad \text{and}\quad C=\max\{||c_1||_\infty, \,||c_2||_\infty \},
\end{equation*}
where $||.||_{\infty}$ stands for the usual $L^\infty$ norm.
In \cite{ASW} and \cite{DVBW}, the authors proposed and studied  a similar linear age-structured model that describes the dynamics of a 
population which consists of proliferating and quiescent cells. It was shown that under some conditions on the support of the transition 
rates between the compartments, the semigroup associated with the linear problem has the property of asynchronous exponential growth 
\cite{ASW,DVBW}. In contrast, in our model individuals are structured by size,  individual development depends on size. 
Another difference is that we consider a very general type of recruitment term (see e.g.~\cite{CS}), namely we assume that 
individuals may have different sizes at birth. We also refer the interested reader to \cite{DHT} and \cite{Heij86} where one-dimensional 
linear size-structured cell population models, with different recruitment terms, were investigated in the framework of positive operator theory. 

\begin{figure}[th]
\begin{center}
\includegraphics[width=70mm]{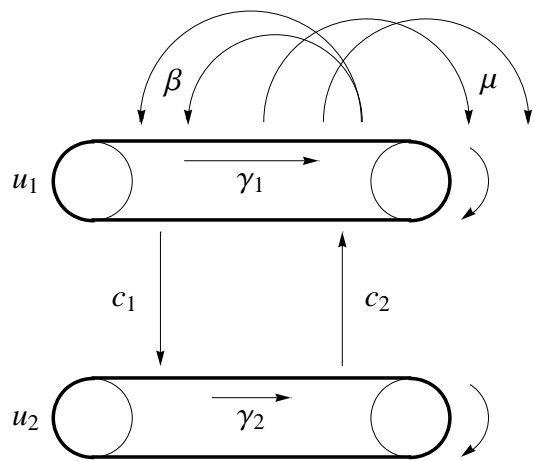}
\caption{The conveyor belts: Schematic representation  of the two-phase population model}\label{fig_1}
\end{center}
\end{figure}

\sector{Existence and positivity of solutions}

Our main objective in this section is to show that our model \eqref{model1}-\eqref{bd_in} is governed by a strongly continuous ($C_0$ for short) positive quasicontractive semigroup on 
the biologically relevant state space $\mathcal{X}=L^1(0,m)\times L^1(0,m)$. This will imply that problem \eqref{model1}-\eqref{bd_in} is well-posed and every solution starting with a non-negative initial condition remains non-negative. 

We equip the space $\mathcal{X}$ with the norm $||\mathbf{x}||_{\mathcal{X}}=||x_1||_1+||x_2||_1$, where $\mathbf{x}=(x_1,x_2)$ and $||\,\cdot\,||_1$ stands for the usual $L^1$ norm. 
For $\mathbf{x}=(x_1,x_2),\,\mathbf{y}=(y_1,y_2)\in\mathcal{X}$ we define the partial ordering $\mathbf{x}\le\mathbf{y}$ if and only if $x_1(s)\le y_1(s)$ and $x_2(s)\le y_2(s)$ for a.e. $s\in [0,m]$. Then $\mathcal{X}$ is a Banach lattice. We refer to the book by Arendt \textit{et al.} \cite[part C]{AGG} for definitions and basic facts about Banach lattices. The dual space of $\mathcal{X}$ is $\mathcal{X}^*=L^\infty(0,m)\times L^\infty(0,m)$ equipped with the norm $||\mathbf{x}^*||_{\mathcal{X}^*}=\sup\{||x_1^*||_\infty,||x_2^*||_\infty\}$, where $\mathbf{x}^*=(x_1^*,x_2^*)$ and $||\,\cdot\,||_\infty$ stands for the usual $L^\infty$ norm. Finally, the natural pairing (or semi-inner product) between elements $\mathbf{x}=(x_1,x_2)\in\mathcal{X}$ and $\mathbf{x^*}=(x^*_1,x^*_2)\in\mathcal{X}^*$ is given by
\begin{equation*}
    \langle \mathbf{x},\mathbf{x^*} \rangle_{\mathcal{X}} = \int_0^m (x_1(s)x^*_1(s) + x_2(s)x^*_2(s))\,\ud s.
\end{equation*}
Next we cast system \eqref{model1}-\eqref{bd_in} in the form of an abstract Cauchy problem on the state space $\mathcal{X}$ as follows 
\begin{equation}\label{abstr}
    \frac{d}{dt}\, \mathbf{u} = \left(\A + \B\right)\,\mathbf{u},\quad    \mathbf{u}(0)=\mathbf{u}_0,
\end{equation}
where $\mathbf{u}=(u_1,u_2)$. The unbounded part is given by 
\begin{equation*}
\A\,\mathbf{u} =  \begin{pmatrix}
- \frac{d}{ds} \left(\gamma_1u_1\right)  \\   -\frac{d}{ds}\left(\gamma_2 u_2\right) \,
   \end{pmatrix}
\end{equation*}
with dense domain
\begin{equation*}
\text{Dom}(\A)=\left\{\mathbf{u}\in W^{1,1}(0,m)\times W^{1,1}(0,m)\,|\,\mathbf{u}(0)=\mathbf{0}\right\}.
\end{equation*}
Notice that here $\mathbf{0}$ is the zero vector in $\RR^2$. The bounded part is given by
\begin{equation*}
\B\,\mathbf{u} =  \begin{pmatrix}
-\left(\mu + c_1\right)u_1 +\int_0^m \beta(\,\cdot\,,y)u_1(y)\,\ud y +c_2u_2 \\  -c_2u_2+c_1u_1
   \end{pmatrix}\quad \text{on}\ \mathcal{X}.
\end{equation*}
Our aim is to establish that the linear operator $\mathcal{A+B}$ is a generator of a positive quasicontractive semigroup (see 
\cite{F,FH4}). To this end first we recall (see e.g. \cite{AGG,CH,NAG}) some basic concepts and results from the theory of linear operators acting on (non-reflexive) Banach spaces. Let $\mathcal{L}$ be a linear operator defined on the real Banach lattice $\mathcal{Y}$ with norm $||.||_{\mathcal{Y}}$. $\mathcal{L}$ is called \textit{dissipative} if for every $\lambda>0$ and $y\in \text{Dom}(\mathcal{L})$, 
\begin{equation*}
 ||(\mathcal{I}-\lambda\mathcal{L})y||_{\mathcal{Y}}\ge ||y||_{\mathcal{Y}}.
\end{equation*} 
A function $f:\mathcal{Y}\to\RR$ is called \textit{sublinear} if 
\begin{align*}
& f(y+z)\le f(y)+f(z),\quad y,z\in\mathcal{Y}\\
& f(\lambda y)=\lambda f(y),\quad \lambda\ge 0,\quad y\in\mathcal{Y}.
\end{align*}
If in addition $f(y)+f(-y)>0$ holds true for $y\ne 0$ then $f$ is called a \textit{half-norm} on $\mathcal{Y}$. The linear operator $\mathcal{L}$ is called $f$-\textit{dissipative} if
\begin{equation*}
f(y)\le f(y-\lambda \mathcal{L}y),\quad \lambda\ge 0,\quad y\in \text{Dom}(\mathcal{L}).
\end{equation*}
An operator $\mathcal{L}$ is called \textit{dispersive}, if it is $p$-dissipative with respect to the canonical half-norm 
\begin{equation}\label{halfnorm}
p(y)=||y^+||_{\mathcal{Y}},
\end{equation}
where $y^+=y\vee 0$ (and $y^-=(-y)^+$). In our setting, the positive part of an element 
$\mathbf{x}=(x_1,x_2)$ of the state space $\mathcal{X}$ is defined as follows 
\begin{align}
\mathbf{x}^+=\left\{ \begin{array}{lll}  (x_1(s),x_2(s)) & \text{if}   & x_1(s)>0,\  x_2(s)>0, \\ 
        (x_1(s),0)  & \text{if}   & x_1(s)>0,\  x_2(s)\le 0, \\
        (0 ,x_2(s)) & \text{if}   & x_1(s)\le 0,\  x_2(s)>0, \\
        (0,0)  & \text{if} & x_1(s)\le 0,\ x_2(s)\le 0 \end{array}
\right. .\label{postivepart}
\end{align}
Clearly, $p$ defined by \eqref{halfnorm} is a continuous sublinear functional and its subdifferential is given by
\begin{equation*}
\ud p\,(y)=\left\{\phi_y\in \mathcal{Y}^*_+\ \text{such that}\ ||\phi_y||_{\mathcal{Y}^*}\le 1,\ \langle y,\phi_y\rangle_{\mathcal{X}}=||y^+||_{\mathcal{Y}}\right\},
\end{equation*}
where $\mathcal{Y}^*_+$ is the positive cone of $\mathcal{Y}^*$. We also note that it follows from the Hahn-Banach Theorem that $\ud p(y)\ne \emptyset$ for every $y\in\mathcal{Y}$. In fact, $\phi_y\in\ud p\,(y)$ if and only if $\phi_y(s)=1$ if $y(s)>0$, $0\le\phi_y(s)\le 1$ if $y(s)=0$ and $\phi_y(s)=0$ if $y(s)<0$.

We recall that a $C_0$ semigroup $\T(t)$ is called \textit{quasicontractive} if 
\begin{equation*}
||\T(t)||\le e^{\omega t},\quad t\ge 0,
\end{equation*}
for some $\omega\in\RR$, and it is called \textit{contractive} if $\omega\le 0$. Quasicontractive semigroups on $L^1$ spaces are of special interest, for example it can be shown (see e.g. \cite{YK}) that every quasicontractive semigroup on an $L^1$ space admits a minimal dominating positive semigroup (so called modulus semigroup, see e.g.~\cite[p. 278]{AGG} for a definition), that itself is  quasicontractive. Let us first recall the following characterization theorem from \cite[Corollary 7.15]{CH} (see also Theorem 1.2 in Sect. C-II in \cite{AGG}).

\begin{theorem}\label{dispersive}
Let $\mathcal{Y}$ be a Banach lattice and let $\mathcal{L}\,:\,\text{Dom}(\mathcal{L})\to\mathcal{Y}$ be a linear operator. The following statements are equivalent.
\begin{enumerate}
\item[(i)] $\mathcal{L}$ is the generator of a positive contractive semigroup. 
\item[(ii)] $\mathcal{L}$ is densely defined, $\text{Rg}(\lambda\mathcal{I}-\mathcal{L})=\mathcal{Y}$ for some $\lambda>0$, and $\mathcal{L}$ is dispersive.
\end{enumerate}
\end{theorem}
We also recall (see e.g. \cite{AGG}) that $\mathcal{L}$ is dispersive if for every $y\in\text{Dom}(\mathcal{L})$ there exists 
$\phi_{y}\in\mathcal{Y}^*$ with $0\le\phi_y$,  $||\phi_y||_{\mathcal{Y}^*}\le 1$ and $\langle y,\phi_y\rangle_{\mathcal{X}}=||y^+||_{\mathcal{Y}}$ such that 
\begin{equation*}
\langle\mathcal{L}y,\phi_y\rangle_{\mathcal{X}}\le 0. 
\end{equation*}
In fact, we have $\phi_y\in\ud p\,(y)$.  Let us finally recall the following generation theorem from \cite[Corollary 5.8]{NAG}.
\begin{theorem}\label{Trotter}
Let $\mathcal{T}(t)$ and $\mathcal{S}(t)$ be strongly continuous semigroups on $\mathcal{X}$ with generators $\mathcal{A}$ and $\mathcal{B}$, respectively, satisfying
\begin{equation*}
\left|\left|\left[\mathcal{T}\left(\frac{t}{n}\right)\mathcal{S}\left(\frac{t}{n}\right)\right]^n\right|\right|\le M\, e^{\omega t},\quad \text{for all}\quad t\ge 0,\,n\in\mathbb{N},
\end{equation*}
for some constants $M\ge 1$ and $\omega\in\mathbb{R}$. Let $D=\text{Dom}(\mathcal{A})\cap\text{Dom}(\mathcal{B})$ and assume that $D$ and $(\lambda_0-\mathcal{A}-\mathcal{B})D$ are dense in $\mathcal{X}$ for some $\lambda_0>\omega$. Then $\mathcal{C}=\overline{\mathcal{A}+\mathcal{B}}$ generates a strongly continuous semigroup $\mathcal{U}(t)$ given by the {\it Trotter product formula}
\begin{equation}\label{productformula}
\mathcal{U}(t)x=\lim_{n\to\infty}\left[\mathcal{T}\left(\frac{t}{n}\right)\mathcal{S}\left(\frac{t}{n}\right)\right]^nx,\quad x\in\mathcal{X},
\end{equation}
with uniform convergence for $t$ in compact intervals.
\end{theorem}
It follows from the previous result that if $\mathcal{T}(t)$ and $\mathcal{S}(t)$ are both positive then also $\mathcal{U}(t)$ is positive. It is also clear from \eqref{productformula} that if $||\mathcal{T}(t)||\le\exp\{\omega_1 t\}$ and $||\mathcal{S}(t)||\le\exp\{\omega_2 t\}$ then we have $||\mathcal{U}(t)||\le \exp\{(\omega_1+\omega_2)t\}$.

Our main result is the following.
\begin{theorem}\label{posgen}
The operator ${\mathcal A} + {\mathcal B}$ generates a positive strongly continuous  quasi\-contractive semigroup $\{{\mathcal T}(t)\}_{t\geq 0}$ of bounded linear operators on ${\mathcal X}$.
\end{theorem}
\Pf{.} 
Our aim is to apply first Theorem \ref{dispersive} to show that the operator $\mathcal{A}$ generates a positive contractive semigroup. To this end, for every $\mathbf{u}\in\text{Dom}(\mathcal{A})$ we define $\mathbf{\phi}_\mathbf{u}\in\mathcal{X}^*$ by 
\begin{align}
\phi_u(s)=\left\{ \begin{array}{lll}  (1,1) & \text{if}   & u_1(s)>0,\  u_2(s)>0, \\ 
        (1,0)  & \text{if}   & u_1(s)>0,\  u_2(s)\le 0, \\
        (0 ,1) & \text{if}   & u_1(s)\le 0,\  u_2(s)>0, \\
        (0,0)  & \text{if} & u_1(s)\le 0,\ u_2(s)\le 0 \end{array}
\right. .\label{pospart}
\end{align}
Then 
\begin{equation*}
||\mathbf{\phi_u}||_{\mathcal{X}^*}\le 1,
\end{equation*}
and clearly 
\begin{equation*}
\langle \mathbf{u},{\bf \phi_u}\rangle_{\mathcal{X}}=\int_0^m \left(u_1(s)\chi_{u_1^+}(s)+u_2(s)\chi_{u_2^+}(s)\right)\,\ud s=||\mathbf{u}^+||_{\mathcal{X}},
\end{equation*} 
where $\chi_{u_1^+}$ and $\chi_{u_2^+}$ stand for the characteristic function of the support of $u_1^+$ and $u_2^+$, respectively.
Making use of \eqref{pospart} we obtain 
\begin{align}
\langle\mathcal{A}\mathbf{u},\mathbf{\phi_u}\rangle_{\mathcal{X}} 
& = -\int_0^m\left(\gamma_1(s)u_1(s)\right)_s\chi_{u_1^+}(s)\,\ud s -\int_0^m\left(\gamma_2(s)u_2(s)\right)_s\chi_{u_2^+}(s)\,\ud s \nonumber \\
& \le -\left(\gamma_1(m)u_1(m)\right)\chi_{u_1^+}(m)-\left(\gamma_2(m)u_2(m)\right)\chi_{u_2^+}(m) \nonumber \\
& \le 0. \label{est3}
\end{align}
Hence the operator $\mathcal{A}$ is dispersive. $\mathcal{A}$ is also densely defined. 
It remains to show the range condition, i.e.~$\text{Rg}(\lambda\I-\A) = \X$.  We observe that the equation
\begin{equation}\label{resequation}
(\lambda \I-\mathcal{A})\,\mathbf{u}=\mathbf{h}
\end{equation}
for $(h_1,h_2)=\mathbf{h}\in\mathcal{X}$ and $\lambda>0$ sufficiently large has a unique solution $(u_1,u_2)=\mathbf{u}\in \text{Dom}({\mathcal A})$, given by
\begin{align} \begin{pmatrix} u_1(s) \\ u_2(s) \end{pmatrix}= \begin{pmatrix} \exp\left\{-\int_0^s\frac{\lambda}{\gamma_1(y)}\,\ud y\right\}\int_0^s\exp\left\{\int_0^y\frac{\lambda}{\gamma_1(z)}\,\ud z\right\}\frac{h_1(y)}{\gamma_1(y)}\,\ud y \\ 
\exp\left\{-\int_0^s\frac{\lambda}{\gamma_2(y)}\,\ud y\right\}\int_0^s\exp\left\{\int_0^y\frac{\lambda}{\gamma_2(z)}\,\ud z\right\}\frac{h_2(y)}{\gamma_2(y)}\,\ud y  
\end{pmatrix}.\label{udef1}
\end{align}
The fact that $\mathbf{u}$  defined by \eqref{udef1} is an element of  $\text{Dom}({\mathcal A})$ follows from
\begin{align*}
|u_i'(s)| & \le \left|\frac{h_i(s)}{\gamma_i(s)}\right|+\frac{\lambda}{\gamma_i(s)}\int_0^m\exp\left\{-\int_y^s\frac{\lambda}{\gamma_i(z)}\,\ud z\right\}\frac{|h_i(y)|}{\gamma_i(y)}\,\ud y \nonumber \\
& \le \left|\frac{h_i(s)}{\gamma_i(s)}\right|+M^i_{\lambda},\quad\quad i=1,2,
\end{align*}
for $\lambda$ large enough for some $M^i_{\lambda}<\infty$, $i=1,2$, that is $\mathbf{u}\in W^{1,1}(0,m)\times W^{1,1}(0,m)$, 
and the range condition is satisfied. Theorem \ref{dispersive} gives that 
$\mathcal{A}$ is a generator of a positive contractive semigroup. Since $\mathcal{B}$ is bounded and 
\begin{equation*}
\mathcal{B}+||\mathcal{B}||\,\mathcal{I}\ge 0
\end{equation*}
it follows (see e.g. \cite[Theorem 1.11, C-II]{AGG}) that $\mathcal{B}$ generates a positive quasicontractive semigroup 
and the proof is completed on the grounds of Theorem \ref{Trotter}. 
\eofproof

\begin{remark}
We rewrite the operator $\mathcal{B}$ as 
\begin{equation*}
\mathcal{B}=\mathcal{B}_1+\mathcal{B}_2,
\end{equation*}
where
\begin{equation*}
\B_1\,\mathbf{u} =  \begin{pmatrix}
-\left(\mu + c_1\right)u_1  \\  -c_2u_2
   \end{pmatrix}\quad \text{on}\ \mathcal{X}.
\end{equation*}
and
\begin{equation*}
\B_2\,\mathbf{u} =  \begin{pmatrix}
\int_0^m \beta(\,\cdot\,,y)u_1(y)\,\ud y +c_2u_2 \\  c_1u_1
   \end{pmatrix}\quad \text{on}\ \mathcal{X},
\end{equation*}
Then we have 
\begin{equation*}
||\exp\{t\mathcal{B}_2\}||\le\exp\{t||\mathcal{B}_2||\}\le\exp\{t(mB+C)\},
\end{equation*} 
and 
\begin{equation*}
||\exp\{t\mathcal{B}_1\}||\le\exp\left\{-t\cdot\min\left\{\inf_{s\in[0,m]}\{\mu(s)+c_1(s)\},\inf_{s\in[0,m]}\{c_2(s)\}\right\}\right\}.
\end{equation*}
Therefore if 
\begin{equation*}
mB+C<\min\left\{\inf_{s\in[0,m]}\{\mu(s)+c_1(s)\},\inf_{s\in[0,m]}\{c_2(s)\}\right\}
\end{equation*}
then the growth bound $\omega_0$ of the semigroup is negative, hence it is uniformly exponentially stable (see e.g. \cite{NAG}), i.e. the population dies out.
\end{remark}

\sector{Asymptotic behavior}

Model \eqref{model1}-\eqref{bd_in} is linear hence one naturally expects that solutions either grow or decay exponentially in time (unless they are at an equilibrium). In this scenario certain stability properties of solutions can be efficiently investigated.   
For example for a simple age-structured population one can often show that solutions behave asymptotically as $n(a,t)\approx Ce^{rt}n_*(a)$, a property called asynchronous exponential growth, where $r$ is the Malthusian parameter and $n_*(a)$ is the stable age-profile. We recall (see e.g. \cite{NAG}) that in the framework of linear semigroup theory a strongly continuous semigroup $S=\left\{{\mathcal S}(t)\right\}_{t\geq 0}$ on a Banach space $\mathcal Y$ with generator ${\mathcal O}$ and \textit{spectral bound}
\begin{equation*}
	s\left({\mathcal O}\right)=\sup \left\{\,\text{Re}(\lambda)\::\:\lambda\in\sigma\left({\mathcal O}\right)\,\right\}
\end{equation*}
is said to exhibit \textit{balanced exponential growth} if there exists a projection $\Pi$ on $\mathcal{Y}$ such that 
\begin{equation}\label{beg}
	\lim_{t\to \infty} \|e^{-s\left({\mathcal O}\right)\,t}\,\mathcal{S}(t)-\Pi\|=0.
\end{equation}
The semigroup $S=\left\{{\mathcal S}(t)\right\}_{t\geq 0}$ is said to exhibit \textit{asynchronous exponential growth}  if it 
exhibits balanced exponential growth with a rank one projection $\Pi$. We further refer the reader to the monographs \cite{AGG,CH,NAG} for basic definitions and notions used throughout this section. Note that balanced exponential growth essentially requires that the spectral bound $s(\mathcal{O})$ is a dominant eigenvalue, and that the semigroup $\mathcal{S}$ is essentially compact, i.e. $\omega_{ess}(\mathcal{S})<s(\mathcal{O})$, where $\omega_{ess}(S)$ stands for the essential growth bound of the semigroup 
$\mathcal{S}$ (see e.g. \cite{NAG}).

Our aim in this section is to carry out a spectral analysis of the governing linear semigroup $\T(t)$. 
First we show that solutions of model \eqref{model1}-\eqref{bd_in} exhibit balanced exponential growth (with a finite dimensional projector $\Pi$). 
In other words, model \eqref{model1}-\eqref{bd_in} admits a finite dimensional global attractor. 
Then we show that under some further biologically relevant conditions on the model parameters the governing semigroup is also irreducible
(see below for the definition) therefore establishing that solutions exhibit asynchronous exponential growth,

\begin{proposition}\label{compact}
The spectrum of $\mathcal{A+B}$ can contain only isolated eigenvalues of finite algebraic multiplicity. \end{proposition}
\Pf{.} Our aim is to show that the resolvent operator $R(\lambda,\mathcal{A+B})$ is compact. Since $\mathcal{B}$ is bounded, it is enough to show that $R(\lambda,\mathcal{A})$ is compact.  Since $W^{1,1}(0,m)$ is compactly embedded in $L^1(0,m)$ by the Rellich-Kondrachov theorem \cite[Theorem 6.3, Part I]{Adams}, the claim follows on the ground of \cite[Proposition II.4.25]{NAG}. 
\eofproof 

\begin{sloppypar} Proposition \ref{compact} implies that the essential spectrum of $\mathcal{A+B}$ is empty. 
This would imply that the governing semigroup $\T(t)$ is essentially compact, i.e. \mbox{$\omega_{ess}(\T)<s(\mathcal{A+B})$} holds, 
if we could establish that the spectral mapping theorem for the semigroup and its generator holds true and show for example that the 
point spectrum $\sigma_P(\mathcal{A+B})$ is not empty. We would like to point out that in general the eigenvalue equation\end{sloppypar}  
\begin{equation*}
(\mathcal{A+B}-\lambda\mathcal{I})\mathbf{u}=0
\end{equation*} 
cannot be solved explicitly due to the transfer terms $c_1$, $c_2$ and due to the very general recruitment term. 
Hence one needs to find an indirect proof to show that for the spectral bound $s(\mathcal{A+B})>-\infty$ holds, indeed. 
In the remarkable paper \cite{Heij86}, Heijmans used the Krein-Rutman  theorem to give such a proof for a linear size-structured model. 
Here we present a different approach that uses a lower bound of the birth process by an operator of rank one.
\begin{theorem}\label{eigenvalue}
The generator $\A+\B$ has a non-empty point spectrum.
\end{theorem}
\Pf{.} We rewrite the Cauchy problem \eqref{abstr} as follows 
\begin{equation}\label{abstr2}
    \frac{d}{dt}\, \mathbf{u} = \left(\A + \B_1+\B_3+\B_4\right)\,\mathbf{u},\quad    \mathbf{u}(0)=\mathbf{u}_0,
\end{equation}
where $\B_3+\B_4=\B_2$ and
\begin{equation*}
\B_3\,\mathbf{u} =  \begin{pmatrix}
\int_0^m\beta(\,\cdot\,,y)u_1(y)\,\ud y  \\  0 \, 
   \end{pmatrix}\quad \text{on}\ \mathcal{X}
\end{equation*}
and
\begin{equation*}
\B_4\,\mathbf{u} =  \begin{pmatrix}
c_2u_2 \\ c_1u_1
   \end{pmatrix}\quad \text{on}\ \mathcal{X}.
\end{equation*}
Fix a nontrivial separable kernel $\beta^*$ that satisfies
\begin{equation*}
0\le \beta^*(s,y) = \beta_1(s) \beta_2(y) \le\beta(s,y)
\end{equation*}
and denote by $\B_3^*$ the corresponding operator with birth process defined by $\beta^*$ (such a $\beta^*$ exists, since $\textnormal{supp}\,\beta$ is not emapty). The eigenvalue equation 
\begin{equation*}
(\A+\B_1+\B_3^*-\lambda\mathcal{I})\mathbf{u}=\mathbf{0}
\end{equation*}
admits a non-trivial solution vector $\mathbf{u}\ne\mathbf{0}$ if and only if $\lambda\in\mathbb{C}$ satisfies the following characteristic equation 
\begin{equation}\label{chareq}
1=K(\lambda)=\int_0^m\beta_2(s)\int_0^s\frac{\beta_1(y)}{\gamma_1(y)}\exp\left\{-\int_y^s \frac{\lambda+\mu(r)+c_1(r)+\gamma_1'(r)}{\gamma_1(r)}\,\ud r\right\}\,\ud y\,\ud s.
\end{equation}
It is easily shown that equation \eqref{chareq} admits a unique real solution $\lambda_*$, 
which is in fact a strictly dominant eigenvalue of $\A+\B_1+\B_3^*$. $\B_1$ is linear and bounded and it follows from \cite[Theorem 1.11 C-II]{AGG} that it generates a positive semigroup. Then it follows from Theorem \ref{Trotter} that $\A+\B_1+\B_3^*$ generates a positive semigroup, since $\B_3^*$ is a positive operator. Clearly, $\B_3-\B_3^*$ is a positive operator and we obtain by Corollary VI.1.11 in \cite{NAG}
\begin{equation}\label{first_lower_bound}
-\infty<s(\A+\B_1+\B_3^*)\le s(\A+\B_1+\B_3).
\end{equation}
Finally, since $\B_4$ is positive and $\A+\B_1+\B_3$ is a generator of a positive semigroup, we have by the same corollary
\begin{equation}\label{increase_spectral_bound}
-\infty<s(\A+\B_1+\B_3)\le s(\A+\B_1+\B_3+\B_4)=s(\A+\B),
\end{equation}
and the proof of Theorem \ref{eigenvalue} is complete.  \eofproof 

\begin{remark}
We note that  $K(0)$ is the net reproduction rate of the population in the active phase when the outflow term 
$c_1$ is interpreted as extra mortality (see e.g. \cite{F}) and when the birth rate $\beta$ is separable. 
Hence if $K(0)>1$ then this dominant eigenvalue $\lambda_*$ is positive, 
while $K(0)<1$ implies that $\lambda_*$ is negative, finally $\lambda_*=0$ if and only if $K(0)=1$, 
that is, the net reproduction rate of individuals in the active phase equals one.  
\end{remark}

\begin{remark}
While the spectral bound $s(\A+\B)$ cannot be obtained explicitly,  it can be determined approximately 
by bounding the kernel $\beta$ (either from above or below)  by finite sums of separable kernels
\begin{equation*}
\beta_n^*(s,y) =\sum_{k=1}^n \beta_{1,k}(s) \beta_{2,k}(y) .  
\end{equation*}
Then, instead of \eqref{chareq}, we obtain the characteristic equation as the determinant of an $n$ by $n$ matrix and 
the leading eigenvalue can be determined at least numerically in case of concrete model ingredients.
\end{remark} 

\begin{lemma}\label{evcompact}
The semigroup $\{{\mathcal T}(t)\}_{t\geq 0}$ generated by the operator $\A+\B_1+\B_3+\B_4$ is eventually compact. 
\end{lemma}
\Pf{.} In the view of the Fr\'echet-Kolmogorov compactness criterion in $L^1$ (see e.g. \cite[Chapter X.1]{Y}) we conclude from
\begin{align*}
|\B_3\mathbf{u}\,(s)-\B_3\mathbf{u}\,(s')| & \le  \int_0^m|\left(\beta(s,y)-\beta(s',y)\right)|\,|u_1(y)|\,\ud y \\
& \le ||\beta(s,\,\cdot\,)-\beta(s',\,\cdot\,)||_\infty\,||u||_1
\end{align*}
and the continuity of $\beta$ (see assumption \eqref{assumptions}) that the operator $\B_3$ is compact. Hence it suffices to consider the operator $\A+\B_1+\B_4$. We define the function $\gamma$ as follows
\begin{equation*}
\gamma(s)=\min\{\gamma_1(s),\gamma_2(s)\}\quad \text{for}\quad s\in [0,m],
\end{equation*}
then $\gamma$ is clearly bounded and strictly positive because of our assumptions on $\gamma_1,\gamma_2$ in \eqref{assumptions}. It follows that $\frac{1}{\gamma}$ is integrable and 
\begin{equation}\label{time}
\tau(s)=\int_0^s\frac{1}{\gamma(r)}\,\ud r
\end{equation}
is greater or equal than the time needed for an individual to grow from size $0$ to size $s$ irrespective of the transit rates $c_1$ and $c_2$. This implies that the semigroup $\T_0(t)$ generated by the operator $\A+\B_1+\B_4$ is nilpotent, in particular for any $\mathbf{u}\in \mathcal{X}^+$ we have 
\begin{equation*}
\T_0(t)\mathbf{u}=\mathbf{0}\quad \text{for}\quad t\ge\tau(m),
\end{equation*}
and the claim follows.
\eofproof

\begin{theorem}\label{theorem:beg}
The semigroup $\T(t)$ generated by the operator $\A+\B$ exhibits balanced exponential growth.
\end{theorem}
\Pf{.}
We have already established the existence of a spectral gap for the generator $\A+\B$, hence it only remains to show that the semigroup generated by $\A+\B$ is eventually norm continuous. Then the Spectral Mapping Theorem holds true and it  follows (see e.g. Corollary VI.1.13 in \cite{NAG}) that the boundary spectrum 
\begin{equation*}
\sigma_+(\A+\B) = \sigma(\A+\B)\cap (s(\A+\B)+i\RR)
\end{equation*}
equals $s(\A+\B)$ which is a pole of the resolvent with finite algebraic (hence geometric) multiplicity by Proposition \ref{compact}. The eventual compactness of the semigroup $\T(t)$ established in Lemma \ref{evcompact} implies eventual norm continuity, hence the proof of Theorem \ref{theorem:beg} is complete. 
\eofproof

\begin{remark}
Here we would like to briefly present another idea to show that the semigroup generated by $\A+\B_1+\B_3+\B_4$ is essentially compact. 
This approach works even in the case when individual size may be arbitrary large, i.e. for models when $s\in[0,\infty)$, 
when one cannot establish eventual compactness of the semigroup and the Spectral Mapping Theorem. 
Since the operator $\B_3$ is compact (and positive as well) and $\A+\B_1+\B_4$ is a generator of a positive semigroup we have
\begin{equation*}
\omega_{ess}(\A+\B_1+\B_3+\B_4)=\omega_{ess}(\A+\B_1+\B_4)\le\omega_0(\A+\B_1+\B_4)= s(\A+\B_1+\B_4).
\end{equation*}
Therefore one only really needs to show that
\begin{equation*}
s(\A+\B_1+\B_4)<s(\A+\B_1+\B_3+\B_4).
\end{equation*}
But this follows from Theorem \ref{eigenvalue} and by noting that the eigenvalue equation 
\begin{equation*}
(\A+\B_1+\B_4-\lambda\mathcal{I})\mathbf{u}=\mathbf{0},\quad \mathbf{u}(0)=\mathbf{0}
\end{equation*}
does not admit a non-trivial solution for any $\lambda\in\mathbb{C}$. 
\end{remark} 

We conclude the section with the following result.
\begin{theorem}\label{aeg}
Assume that there exists an $\varepsilon_0>0$ such that for all $0<\varepsilon\le\varepsilon_0$
\begin{equation}\label{assume_birth}
 \int_0^\varepsilon \int_{m-\varepsilon}^m  \beta(s,y)\,\ud y\,\ud s > 0
\end{equation}
and that the transition rates satisfy
\begin{equation}\label{assume_c12}
\inf\, \textnormal{supp} \,c_1 = 0, \quad \text{ and } \quad \sup\,\textnormal{supp} \,c_2 = m.
\end{equation}
Then the semigroup $\T(t)$ generated by $\A+\B$ exhibits asynchronous exponential growth.
\end{theorem}
\Pf{.}
It only remains to show that under conditions \eqref{assume_birth}, \eqref{assume_c12} the semigroup $\T(t)$ is \textit{irreducible},  i.e.~for every $\mathbf{0}\ne\mathbf{u}\in\mathcal{X}_+$ and $\mathbf{0}\ne\mathbf{u}^*\in\mathcal{X}^*_+$ there exists a $t_0$ such that 
\begin{equation}\label{irreducibledef}
\langle \T(t_0)\mathbf{u},\mathbf{u}^*\rangle_{\mathcal{X}}>0.
\end{equation}
Let $\mathbf{u}=(u_1,u_2)$ and $\pi_1$ and $\pi_2$ be the projections onto the first and second coordinates, respectively. First assume $u_1(s)>0$ for $s$ in some interval $[s_0^-,s_0^+]$. Since $\gamma_1>0$, there exist $t_*<t^*$ such that 
\begin{equation*}
\textnormal{supp} \,\pi_1\left(\T(t) \mathbf{u}\right) \cap \textnormal{supp}\,\beta(s,\,\cdot\,)\neq \emptyset
\end{equation*}
\begin{sloppypar}
\noindent for every $t_*\le t\le t^*$ and every $s\in(0,\varepsilon]$. By assumption \eqref{assume_birth}, \mbox{$\pi_1(\T(t) \mathbf{u})(s)>0 $} for $t_*\le t\le t^*$ and $s\in(0,\varepsilon]$ and eventually in a neighborhood of every $s\in(0,m]$. If $u_1(s)=0$ for all $s$, then $u_2(s)>0$ for $s$ in some interval $[s_0^-,s_0^+]$. Since $\gamma_2>0$, there exist $t_*<t^*$ such that 
\begin{equation*}
\textnormal{supp} \,\pi_2\left(\T(t) \mathbf{u} \right)\cap \textnormal{supp}\,c_2\neq \emptyset
\end{equation*}
for every $t_*\le t\le t^*$ and hence by assumption \eqref{assume_c12} $\pi_1\left(\T(t) \mathbf{u} \right)(s)>0$ for $s\in \textnormal{supp}\,c_2$. 
By a similar argument as above,  $\pi_1\left(\T(t) \mathbf{u} \right)(s)>0$ for $s\in (0,\varepsilon]$ and $t$ sufficiently large. 
Finally,  assumption \eqref{assume_c12} also guarantees that  $\pi_2\left(\T(t) \mathbf{u} \right)(s)>0$ for 
$s\in (0,\varepsilon]$ and $t$ sufficiently large.
\end{sloppypar}
\eofproof

\sector{Concluding remarks}
Quiescent phases have been introduced in a variety of biological models, see e.g.~\cite{ASW,DVBW,Had08a,Had08b} and the references therein. 
They arise in biological populations that show distinguished classes of behaviors. For example, it is well known that cells in a tissue 
can either migrate or divide, but not both. To the best of our knowledge, a size-structured models with a quiescent phase has not been studied yet. 

Asynchronous exponential growth is often observed in linear models of population dynamics. 
It represents the fact that although a population is still exponentially growing (or decaying), it's age- or size distribution is ``at equilibrium''. 
Several authors have proved the property of asynchronous exponential growth for age-structured  populations \cite{ASW,DVBW}, 
with a possibly delayed birth process \cite{PT}. In this paper, using the techniques established for semigroups of positive operators, 
we have proved a corresponding result for  a size-structured population with a quiescent compartment. It should be stated here, 
that the conditions that distinguish asynchronous exponential growth from merely balanced  exponential growth, 
namely \eqref{assume_birth} and \eqref{assume_c12} are natural in the sense that they can be achieved by an appropriate choice of the size space. 
If, for example, the possibility to return to the active class would cease at some $m_1<m$, i.e.~$c_2(s)=0$ for all $s\ge m_1$, 
it would be possible to reduce the size space in the second coordinate to $[0,m_1]$, without affecting the asymptotic behavior.

It was shown recently by Hadeler and Thieme  \cite{Had08c} for finite-dimensional models that coupling to a quiescent phase can 
shift spectral bound of a matrix (i.e.~the intrinsic growth rate) in both directions, depending on the transition rates between 
active and quiescent phase. It can be seen easily that the same holds for our infinite-dimensional setting. 
The transition to the quiescent phase can either circumvent a region of high mortality or high reproductive activity in the size 
space of the active population. 

We would like to note that the next natural step will be to incorporate some interaction variables into our model. 
Interaction between individuals in the population may be induced for example by scramble competition for available resources.  
It also seems reasonable for example to require that the transition rates $c_1$ and $c_2$ depend on the standing population sizes 
in the active and resting phases, respectively. Then our model becomes a nonlinear one, and the question arises what types of nonlinearities 
will preserve the asymptotic properties of solutions discussed in the previous section. We note that results in the literature 
for nonlinear models are rather rare, see for example \cite{GW}.

\subsection*{Acknowledgments}

JZF acknowledges financial support from the IMA while being a longterm visitor during 01/2009-02/2009 at the University of Minnesota. 
JZF is also thankful to the Centre de Recerca Mathem\`{a}tica, Universitat Aut\`{o}noma de Barcelona for their hospitality while being a 
participant in the research programme "Mathematical Biology: Modelling and Differential Equations" during 01/2009-06/2009.
PH was supported by IMA postdoctoral fellowship. We thank the referee for valuable comments and suggestions which lead to an improved version of the paper.

\end{document}